\documentclass[letterpaper,11pt]{article}
\usepackage{amsmath,amsthm,amssymb}
\usepackage{graphicx}

\begin{document}

\setlength{\intextsep}{\floatsep}
\def\ldotsplus{\mathinner{\ldotp\ldotp\ldotp\ldotp}}
\setlength{\intextsep}{\floatsep}
\def\ldotscomm{\mathinner{\ldotp\ldotp\ldotp,}}
\def\fourdots{\relax\ifmmode\ldotsplus\else$\m@th \ldotsplus\,$\fi}
\def\dotcoms{\relax\ifmmode\ldotscomm\else$\m@th \ldotsplus\,$\fi}

\long\def\symbolfootnote[#1]#2{\begingroup%
\def\thefootnote{\fnsymbol{footnote}}\footnote[#1]{#2}\endgroup}

\theoremstyle{plain}
\newtheorem{theorem}{Theorem}
\newtheorem{lemma}[theorem]{Lemma}
\newtheorem{corollary}[theorem]{Corollary}
\theoremstyle{definition}
\newtheorem{exmp}[theorem]{Example}
\newtheorem{remark}[theorem]{Remark}

\newcommand{\rarrow}{\overrightarrow}
\newcommand{\ceilrs}{\lceil r/s \rceil}
\newcommand{\fracrs}{\frac{r}{s}}
\newcommand{\fracab}{\frac{a}{b}}

\newcommand{\tpqrs}{T_{p,q}(r,s)}
\newcommand{\tgen}[4]{T_{#1,#2}(#3,#4)}

\newcommand{\spri}{s'}
\newcommand{\rpri}{r'}

\newcommand{\circu}{\operatorname{circ}}
\newcommand{\per}{\operatorname{per}}
\newcommand{\Z}{\mathbb{Z}}
\newcommand{\sgn}{\mbox{sgn}}

\newcommand{\myfig}[2]{\begin{figure}[htbp]\begin{center}
      {\scalebox{.4}{\includegraphics{#1}}}
      \caption{#2}\label{fig:#1}
    \end{center}\end{figure}}

\title{The Combinatorics of a Three-Line Circulant Determinant}
\date{}
\maketitle
\begin{center}
{\large Nicholas A. Loehr,\symbolfootnote[1]{Supported by NSF Postdoctoral
research grants}
 Gregory S. Warrington$^*$ and  Herbert S. Wilf}\\

\vspace{.4in}

{Department of Mathematics, University of Pennsylvania\\ Philadelphia, PA
19104-6395}\\

\vspace{.3in}

{\small\texttt{<nloehr@math.upenn.edu> <gwar@math.upenn.edu>
<wilf@math.upenn.edu>}}
\end{center}

\vspace{.5in}
\begin{abstract}
  We study the polynomial
  $\Phi(x,y)=\prod_{j=0}^{p-1}\left(1-x\omega^j-y\omega^{q j}\right)$,
  where $\omega$ is a primitive $p$th root of unity.  This polynomial
  arises in CR geometry \cite{dan}.  We show that it is the
  determinant of the $p\times p$ circulant matrix whose first row is
  $(1,-x,0,\ldots,0,-y,0,\dotcoms 0)$, the $-y$ being in position $q+1$.
  Therefore, the coefficients of this polynomial $\Phi$ are integers
  that count certain classes of permutations. We show that all of the
  permutations that contribute to a fixed monomial $x^ry^s$ in $\Phi$
  have the same sign, and we determine that sign. We prove that a
  monomial $x^ry^s$ appears in $\Phi$ if and only if $p$ divides
  $r+sq$.  Finally, we show that the size of the largest coefficient
  of the monomials in $\Phi$ grows exponentially with $p$, by proving
  that the permanent of the circulant whose first row is
  $(1,1,0,\dotcoms 0,1,0,\dotcoms 0)$ is the sum of the absolute
  values of the monomials in the polynomial $\Phi$.
\end{abstract}

\newpage

\section{Introduction and statement of results}
The stimulus for this work lies in the study \cite{dan} by John
D'Angelo of invariant holomorphic mappings on hypersurfaces. In that
work a construction is given of a certain real-analytic function
$\Phi$ from which one can define the desired invariant mappings. As a
source of examples the author used the familiar lens spaces ${\cal
  L}(p,q)$, and he showed that the invariant function in \cite{dan}
determines a polynomial in two real variables we call $\Phi$.
Specifically,
\begin{equation}\label{eq:fdef}
  \Phi(x,y)=\Phi_{p,q}(x,y)=\prod_{j=0}^{p-1}\left(1-x\omega^j-y\omega^{q
j}\right),
\end{equation}
where $\omega$ is a primitive $p$th root of unity. For example,
\[\Phi_{8,3}(x,y)=1 - x^8 - 8\,x^5\,y - 12\,x^2\,y^2 + 2\,x^4\,y^4 -
  8\,x\,y^5 - y^8.\]
  
  Hence in the case of lens spaces, $\Phi$ is a polynomial in $x,y$
  that has certain interesting extremal properties. For further
  investigation it is desirable to know more about these polynomials.
  In particular,
\begin{enumerate}
\item \label{it:one}Are its coefficients always integers?
\item \label{it:two}If so, what integers are they?
\item \label{it:three}Precisely which monomials in $x,y$ appear in
$\Phi_{p,q}(x,y)$?
\item \label{it:four}Which of the coefficients of the monomials that appear
are positive and which are negative?
\end{enumerate}

Question \ref{it:one} was already answered in the affirmative in
\cite{dan}.  In Section \ref{sec:circ} we will give a particularly
simple proof (and a combinatorial interpretation to the coefficients),
by exhibiting $\Phi_{p,q}$ as the determinant of a certain $p\times p$
matrix that has integer entries.

Question \ref{it:two} is harder. As a partial answer, in Section
\ref{sec:uniq} as a corollary of Lemma \ref{lem:sign}, we will prove
the following:
\begin{theorem}
\label{th:one}
 In the expansion of the polynomial
\[\Phi_{p,q}(x,y)=\sum_{r,s}a_{p,q}(r,s)x^ry^s\]
the coefficient $a_{p,q}(r,s)$ is equal, aside from its sign, to the
number of permutations $\sigma$ of $p$ letters such that the
differences
\[\{(\sigma(j)-j)\ \mathrm{mod}\, p\}_{j=1}^p\]
take the values $0$, $1$, and $q$ with respective multiplicities
$p-r-s$, $r$, and $s$. Furthermore, these permutations all have the
same signs, and in fact, all have the same cycle type.
\end{theorem}

Regarding Question \ref{it:three}, we obtain the following from Lemma
\ref{lem:nonem} of Section \ref{sec:circ} and Theorem
\ref{thm:nonempty} of Section \ref{sec:exis}:
\begin{theorem}
\label{th:two}
The monomials $x^ry^s$ that appear in $\Phi_{p,q}(x,y)$ (i.e., that
have nonzero coefficients) are precisely those for which $p$ divides
$r+sq$.
\end{theorem}
That $p$ must divide $r + sq$ for $x^ry^s$ to appear with nonzero
coefficient is by far the easier implication to prove.  This necessity
follows from the underlying geometry (see \cite{dan}) or, as we will
show, from a simple counting argument.

Finally, Question \ref{it:four} about the signs of the terms is
settled by the following result which follows from Lemma
\ref{lem:sign} in Section \ref{sec:uniq}.
\begin{theorem}
\label{th:three}
Let $a_{p,q}(r,s)x^ry^s$ be a monomial that appears in $\Phi_{p,q}(x,y)$.
Then the sign of this monomial is positive (resp. negative) if the integer
\[\mathrm{gcd}\left(r,s,\frac{r+sq}{p}\right)\]
 is even (resp. odd).
\end{theorem}

Finally in Section \ref{sec:maxc} we show that, for fixed $q$, the
coefficients in $\Phi_{p,q}$ grow exponentially with $p$.

\begin{remark}
  D'Angelo \cite{danpre} shows that the polynomial $f(x,y)=1 -
  \Phi(x,y)$ is congruent to $(x+y)^p \pmod{p}$ if and only if $p$ is
  prime.
\end{remark}

\begin{remark}
One can also consider expressions of the form
\begin{equation}
  \Theta_{p,q,t} = \prod_{j=0}^{p-1} \left(1-x\omega^{tj}-y\omega^{qj}\right).
\end{equation}
These can be realized as determinants of $p\times p$ matrices of the
form
\begin{equation}\label{eq:geneq}
  \circu(1,0,\dotcoms 0,-x,0,\dotcoms 0,-y,0,\dotcoms 0)
\end{equation}
where the $-x$ and $-y$ appear in the $(t+1)^{\mathrm{st}}$ and
$(q+1)^{\mathrm{st}}$ positions, respectively.  If $\omega^t$ is a
primitive root of unity (i.e., $\gcd(t,p)=1$), then $\omega^q =
\omega^{tq'}$ for some $q'$.  This implies that $\Theta_{p,q,t}$
equals $\Phi_{p,q'}$.  (A similar statement can be made when
$\gcd(q,p) = 1$.)  This extends somewhat the set of $(p,q,t)$ to which
our results apply, but the general case remains open.

The permanents of the $(0,1)$-matrices associated to the
$\Theta_{p,q,t}$ are investigated in \cite{ccr} (see, in particular,
Lemma 12).  We note that, according to Theorem \ref{th:one} above, all
permutations that contribute to a given monomial have the same sign.
Since there is no cancellation, we obtain the following:

\begin{corollary}
  The \textit{permanent} of a $p\times p$ circulant matrix whose first
  row has 1's in columns 1,2, and $q+1$ (and 0's elsewhere) is equal
  to the sum of the absolute values of the coefficients of the
  monomials that occur in $\Phi$.
\end{corollary}
\end{remark}

\section{Circulant matrices}
\label{sec:circ}
A $p\times p$ \textit{circulant matrix} is a matrix of the form
\[C=\left[
\begin{array}{ccccc}
a_0&a_1&a_2&\dots&a_{p-1}\\
a_{p-1}&a_0&a_1&\dots&a_{p-2}\\
a_{p-2}&a_{p-1}&a_0&\dots&a_{p-3}\\
\vdots&\vdots&\vdots&\dots&\vdots\\
a_1&a_2&a_3&\dots&a_0
\end{array}
\right].\]
Since such a matrix is completely specified by, for example, its first row,
we will sometimes refer to it as $\circu(a_0,a_1,\dots ,a_{p-1})$.
A circulant matrix can be written as $C=g(C_0)$ where
$C_0=\circu(0,1,0,\dots,0)$ and
$g(t)=a_0+a_1t+a_2t^2+\dots+a_{p-1}t^{p-1}.$ Since the eigenvalues of
$C_0$ are the $p$th roots of unity, the eigenvalues of the
general circulant matrix $C$ are $g(\omega)$, where $\omega$ runs
through the $p$th roots of unity. Consequently the determinant of any
circulant matrix is the product of these eigenvalues, namely
\[\det{C}=\prod_{\omega^p=1}g(\omega).\]
The above observations are from well known, classical theory of circulant
matrices. See, for example \cite{da}.

If we take $g(t)=1-xt-yt^q$ we see that the polynomial $\Phi(x,y)$, whose
study is the main object of this paper,
is the determinant of $g(C_0)$, as stated above. {F}rom the form of $g$ we
see at once that the polynomial $\Phi$ has integer coefficients, thus
answering Question \ref{it:one} by inspection.

If we write $\Phi(x,y)=\sum_{r,s}a(r,s)x^ry^s$, then we can give a
combinatorial interpretation to the coefficients $a(r,s)$.  Indeed, by
expanding the circulant determinant
\[\Phi(x,y)=\det{(I-xC_0-yC_0^q)}=
\left|
\begin{array}{cccccccc}
1&-x&0&\dots&0&-y&0&0\\
0&1&-x&\dots&0&0&-y&0\\
0&0&1&\dots&0&0&0&-y\\
\vdots&\vdots&\vdots&\dots&\vdots&\vdots&\vdots&\vdots\\
-x&0&0&\dots&-y&0&0&1
\end{array}
\right|,\]

we see that the coefficient of $(-1)^{r+s}x^ry^s$ is the sum of the
signs of those permutations of $p$ letters that ``hit'' $r$ of the
$x$'s in the matrix and $s$ of the $y$'s, the remaining values being
fixed points. Thus, let $\tpqrs$ denote the set of all permutations
$\sigma$ of $1,2,\dots,p$ such that

\begin{enumerate}
\item $\sigma$ has exactly $p-r-s$ fixed points, and
\item for exactly $r$ values of $j$ we have $\sigma(j)-j$ congruent to 1
modulo $p$, and
\item for exactly $s$ values of $j$ we have $\sigma(j)-j$ congruent to $q$
modulo $p$.
\end{enumerate}

Then $(-1)^{r+s}a(r,s)$ is the excess of the number of even
permutations in $\tpqrs$ over the number of odd permutations in
$\tpqrs$.

As an example, take $p=5$ and $q=3$. Then
\[\Phi(x,y)=1-x^5-5x^2y-5xy^3-y^5.\]
Let's check the coefficient of $x^2y$. The set $\tgen{5}{3}{2}{1}$ consists
of the following permutations of $5$ letters:
\[ \{1,2,4,5,3\},\{1,3,4,2,5\},\{2,3,1,4,5\},\{2,5,3,4,1\},\{4,2,3,5,1\}.\]
These are all even permutations, hence $-a(2,1)$ is $5$, as we also see by
inspection of $\Phi$.
Note that all of the permutations in $\tgen{5}{3}{2}{1}$
have the same cycle structure, viz. a 3-cycle and two fixed points.

Our goal is to show the following:
\begin{itemize}
\item (Uniqueness) If $\tpqrs$ is nonempty, then every $\sigma\in
  \tpqrs$ has the same cycle structure.  We will explicitly describe
  this cycle structure.
\item (Existence) $\tpqrs$ is nonempty if and only if $p$ divides
  $r+sq$.
\end{itemize}

We first consider two special cases. If $r=s=0$, then $\tgen{p}{q}{0}{0}$
consists of the identity permutation.  If $s=0$ and $r>0$, it is easy
to see from the definitions that $\tgen{p}{q}{r}{0}$ is nonempty iff $r=p$,
in which case the only element of this set is the cycle
$(1,2,\ldots,p)$.  In what follows, therefore, we assume $s > 0$.

\section{Uniqueness of cycle structure}
\label{sec:uniq}

It is convenient to introduce the following notation for a permutation
$\sigma\in \tpqrs$. Write $\sigma$ uniquely (up to order) as a product
of $k\geq 0$ disjoint cycles $C_1,\ldots,C_k$ of lengths greater than
1.  If $k=0$, then $\sigma$ is the identity. This happens only in the
trivial case $r=s=0$, so we assume $k>0$ from now on.

We will represent each cycle $C_i$ by a pair $(x_i;w_i)$, where
$x_i\in\{1,2,\ldots,p\}$ and $w_i$ is a word consisting
of $r_i$ 1's and $s_i$ $q$'s. Here, $x_i$ is an arbitrary
point appearing in the cycle $C_i$, $r_i+s_i$ is the number
of points involved in the cycle, and the word $w_i$ gives
the differences (mod $p$) between consecutive elements of the
cycle starting at $x_i$.
Formally, if $w_i=w_i(1)w_i(2)\cdots w_i(r_i+s_i)$, then
\begin{equation}\label{eq:cycle-formula}
 \sigma^t(x_i)\equiv x_i+\sum_{j=1}^t w_i(j) \pmod{p},
\mbox{ for $0\leq t\leq r_i+s_i$.}
\end{equation}
(We take our residue system modulo $p$ to be the set
$\{1,2,\dotcoms p\}$.)  For example, when $q=3$ and $p=10$, the pair
$(4; 3,1,1,3,1,1)$ represents the cycle $(4,7,8,9,2,3)$. The pair $(8;
1,3,1,1,3,1)$ also represents this cycle.

\begin{lemma}\label{lem:nonem}
  If $\tpqrs$ is nonempty, then $p$ divides $r+sq$.
\end{lemma}
\begin{proof}
Take any $\sigma\in \tpqrs$, and describe $\sigma$ using
the notation above. Each cycle $C_i$ has $r_i+s_i$ elements in it.
Letting $t=r_i+s_i$ in (\ref{eq:cycle-formula}) gives
\[ x_i=\sigma^{r_i+s_i}(x_i)\equiv x_i+\sum_{j=1}^{r_i+s_i} w_i(j)
 \equiv x_i+r_i\cdot 1+s_i\cdot q \pmod{p}. \]
Thus, $p$ divides $r_i+qs_i$ for each $i$.
It is easy to see from the definitions that $r=r_1+\cdots+r_k$
and $s=s_1+\cdots+s_k$. Hence,
$r+qs=\sum_{i=1}^k (r_i+qs_i)$ is also divisible by $p$.
\end{proof}

By the proof of the last lemma, $p$ divides all the quantities
$r_i+qs_i$. So, given $\sigma\in \tpqrs$,
we can define positive integers $\ell_i=(r_i+qs_i)/p$
and $\ell=(r+qs)/p=\sum_{i=1}^k \ell_i$.

\begin{lemma}\label{lem:gcdone}
  If $\tpqrs$ is nonempty, then $\gcd(r_i,s_i,\ell_i)=1$ for $1\leq
  i\leq k$.
\end{lemma}
\begin{proof}
Fix $i$ between $1$ and $k$.
We assume that $\gcd(r_i,s_i,\ell_i)=d>1$ and derive a contradiction.
Set $\rpri =r_i/d$, $\spri =s_i/d$, and $\ell'=\ell_i/d$.
Since $r_i+qs_i=\ell_i p$, we have $\rpri +q\spri =\ell' p$.

We claim that there exists a string of $\rpri +\spri $ consecutive symbols in
$w_i$ consisting of $\rpri $ 1's and $\spri $ $q$'s.  To prove this, we
start by factoring the word $w_i$ into $d$ subwords
\[ w_i=v_1v_2\cdots v_d, \]
where each word $v_j$ has length $\rpri +\spri $.
For $1\leq j\leq d$, let $v_j$ consist of $a_j$ 1's and $b_j$ $q$'s,
where $a_j+b_j=\rpri +\spri $.  If $a_j=\rpri $ for any $j$, then the claim is true.
If $a_j>\rpri $ for all $j$, then the total number of 1's in $w_i$
is greater than $\rpri d=r_i$, which is a contradiction.
If $a_j<\rpri $ for all $j$, then the total number of 1's in $w_i$
is less than $\rpri d=r_i$, which is a contradiction.
So we are reduced to the case where $a_{j_1}>\rpri $ for some $j_1$
and $a_{j_2}<\rpri $ for some $j_2$. Clearly, in this case we can
choose $j_1$ and $j_2$ with $|j_2-j_1|=1$. We have (say)
\begin{eqnarray*}
 v_{j_1} &=& x_1x_2\cdots x_{\rpri +\spri }, \\
v_{j_2}=v_{j_1+1} &=& x_{\rpri +\spri +1}\cdots x_{2\rpri +2\spri }.
\end{eqnarray*}
Define a function $g:\{1,2,\ldots,\rpri +\spri +1\}\rightarrow\Z$
by declaring $g(m)$ to be the number of 1's in the string
$x_mx_{m+1}\cdots x_{m+\rpri +\spri -1}$. Then $g(1)=a_{j_1}>\rpri $ and
$g(\rpri +\spri +1)=a_{j_2}<\rpri $ and $|g(i+1)-g(i)|\leq 1$ for all $i$.
Hence, there must exist some $m$ with $g(m)=\rpri $. The subword of $w_i$
of length $\rpri +\spri $ beginning with $x_m$ must then contain
$\rpri $ 1's and $\spri $ $q$'s.  This proves the claim.

By the claim, for some $j\geq 0$ there is a subword
\[ w_i(j+1),w_i(j+2),\ldots,w_i(j+\rpri +\spri ) \]
consisting of $\rpri $ 1's and $\spri $ $q$'s. Consider the elements
\[ y=\sigma^j(x_i),\ z=\sigma^{j+\rpri +\spri }(x_i) \]
on the cycle $C_i$. On one hand, we have $y\neq z$
since $\rpri +\spri =(r_i+s_i)/d$ is less than the length $r_i+s_i$ of $C_i$.
On the other hand, (\ref{eq:cycle-formula}) gives
\[ z-y\equiv\sum_{m=j+1}^{j+\rpri +\spri } w_i(m)
   \equiv \rpri +\spri q =\ell' p \equiv 0\pmod{p}. \]
Since $1\leq y,z\leq p$, we get $y=z$, a contradiction.
\end{proof}


We will now precisely characterize the cycles in $C$.  In order
to avoid having to keep track of when $z+q \leq p$ in what follows, we
introduce the following notation: For $z_1,\dotcoms z_m\in [p]$ with
$m \geq 3$, we write $\rarrow{z_1\cdots z_m}$ if there exists a $j$
with $1\leq j \leq m$ such that
\begin{equation}\label{eq:gr2}
  z_j < \cdots < z_m < z_1 < \cdots < z_{j-1}.
\end{equation}
If we think of $[p]$ being arranged in clockwise order around a
circle, then $\rarrow{z_1\cdots z_m}$ amounts to having the clockwise
traversal of $z_1$ to $z_m$ encounter $z_i$ before $z_j$ if and only
if $i < j$.

\begin{lemma}\label{lem:grcyc}
  Let $z_1,\dotcoms z_m\in [p]$ and set $\pi(z) = z + q \pmod{p}$.
  Then
  \begin{equation}\label{eq:gr7}
    \rarrow{z_1\cdots z_m} \Longleftrightarrow
    \rarrow{\pi(z_1)\cdots \pi(z_m)}.
  \end{equation}
\end{lemma}

\begin{proof}
  Assume $\rarrow{z_1\cdots z_m}$ and pick $j$ as in \eqref{eq:gr2}.
  Certainly
  \begin{equation}\label{eq:gr3}
    z_j + q < \cdots < z_m + q < z_1 + q < \cdots < z_{j-1} + q.
  \end{equation}
  If $z_j + q > p$ or $z_{j-1} + q \leq p$, then we immediately obtain
  $\rarrow{\pi(z_1)\cdots \pi(z_m)}$.  Otherwise, there is a
  minimal $t$ (with respect to the order $j < \cdots < m < 1 < \cdots
  < j-1$), $t \neq j$, such that $z_t + q > p$.  Then the only
  nontrivial inequality in
  \begin{equation}\label{eq:gr4}
    \pi(z_t) < \cdots < \pi(z_{j-1}) < \pi(z_j) < \cdots < \pi(z_{t-1})
  \end{equation}
  is $\pi(z_{j-1}) < \pi(z_j)$.  But this must be true as $z_{j-1} - p
  \leq 0 < z_j$ implies $\pi(z_{j-1}) = z_{j-1} + q - p < z_j + q =
  \pi(z_j)$.  The other implication of \eqref{eq:gr7} results from the
  above arguments applied to $\pi^{-1}$, which is the map sending $z$
  to $z+p-q$ (mod $p$).
\end{proof}

\begin{lemma}\label{lem:grinorder}
  For $\sigma\in \tpqrs$, we must have $r_1=r_2=\cdots=r_k$ and
  $s_1=s_2=\cdots=s_k$.
\end{lemma}

\begin{proof}
  Let $C_k$ and $C_l$ be two distinct cycles in $\tpqrs$.  For
  simplicity, we substitute $C,C',a,b,a',b\,'$ for $C_k,C_l,
  r_k,s_k,r_l,s_l$, respectively.  Write
  \begin{equation}\label{eq:grcxv}
    C = (x;v) \text{ where } v = 1^{\alpha_1}q\cdots 1^{\alpha_b}q,
    x\in [p], \text{ and } \sum_i \alpha_i = a.
  \end{equation}
  In traversing the orbit of $x$ under $C$, we will refer to those $z$
  for which $C(z) \equiv z + q \pmod{p}$ as ``$q$-steps''; ``$1$-steps''
  are defined analogously.

  If $b$ were to be $0$, then $a$ would equal $p$ and $C =
  (1,2,\dotcoms p)$.  In this scenario, there are no nontrivial cycles
  disjoint from $C$.  This contradicts our hypothesis.  Hence, $b >
  0$.  Similarly, $b\,' > 0$.  We wish to show that $b\,' \geq b$.  If $b
  = 1$, there is nothing to prove, so assume furthermore that $b > 1$.

  Set $d_1 = x$ and $e_1 = C^{\alpha_1}(x)$.  Then, for $2 \leq i \leq
  b$, we recursively define $d_i = C(e_{i-1})$ and $e_i =
  C^{\alpha_i}(d_i)$.  Note that $d_i$ is the image of the
  $(i-1)^{\mathrm{st}}$ $q$-step, and $e_i$ is the $i$-th $q$-step.

  There exists a unique permutation $\tau$ such that $\tau(1) = 1$ and
  \begin{equation}\label{eq:grperm}
    \rarrow{d_{\tau(1)}e_{\tau(1)}\cdots
      d_{\tau(b)}e_{\tau(b)}}.
  \end{equation}
  Notice that each $e_{\tau(j)}$ is a $q$-step of $C$.  Now let $z$
  be moved by $C'$ (hence fixed by $C$).  For brevity in what follows,
  we interpret the indices of $e$ and $d$, and the arguments of $\tau$,
  modulo $b$.  Set
  \begin{align*}
    V_{\tau(j)}
    &= \{y\in [p]: C(y) = y \text{ and }
    \rarrow{e_{\tau(j)}ye_{\tau(j+1)}}\}\\
    &= \{y\in [p]: C(y) = y \text{ and }
    \rarrow{d_{\tau(j)}yd_{\tau(j+1)}}\}.
  \end{align*}
  The equality of these two sets is due to the fact that each of the cyclic
  intervals $\{z:\rarrow{d_{\tau(j)} z e_{\tau(j)}}\}$ consists only of
points
  moved by $C$.

  By \eqref{eq:grperm}, $z\in V_{\tau(j)}$ for a unique $j$.  If $z$
  is a $1$-step of $C'$, then $C'(z) \in V_{\tau(j)}$ also as $C$
  and $C'$ are disjoint.  If $z$ is a $q$-step of $C'$, then $\pi(z) =
  C'(z)$.  So by Lemma \ref{lem:grcyc}, since
  $\rarrow{e_{\tau(j)}ze_{\tau(j+1)}}$, we find that
  $\rarrow{d_{\tau(j)+1}C'(z)d_{\tau(j+1)+1}}$.  Or, equivalently,
  that $C'(z) \in V_{\tau(j)+1}$.  Iterating this argument, we see
  that the orbit of $z$ visits $V_{\tau(j)}, V_{\tau(j+1)},
  V_{\tau(j+2)}\ldots $ in turn.  We conclude that $C'$ has at least
  $b$ $q$-steps.  Then, by definition, $b\,' \geq b$.  Arguing with the
  roles of $C$ and $C'$ switched, we find that $b = b\,'$.

  To show that $a = a'$, it suffices to consider the equalities $a +
  bq = \ell p$ and $a' + bq = \ell'p$.  Subtracting, $a-a' =
  (\ell-\ell')p$.  Since $b = b\,' > 0$, we know that $0\leq a,a' < p$.
  So $-p < a-a' < p$.  It follows that $a = a'$.
\end{proof}

\myfig{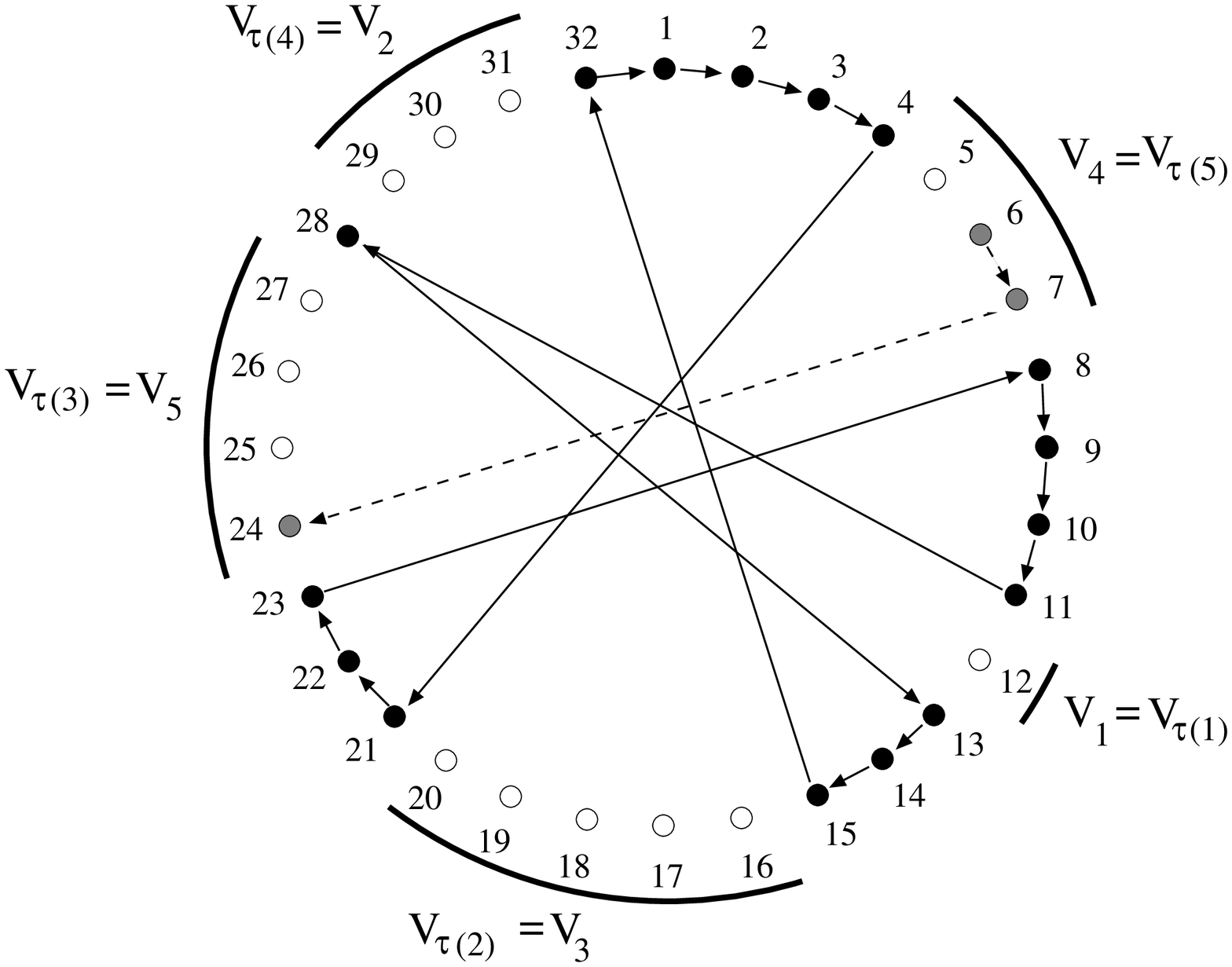}{Illustration for Example \ref{ex:cir}}

\begin{exmp}\label{ex:cir}
  Set $p = 32$ and $q = 17$.  The cycle
  \begin{equation*}
    (8,9,10,11,28,13,14,15,32,1,2,3,4,21,22,23)
  \end{equation*}
  illustrated in Figure \ref{fig:circle} can be written according to
  the conventions of \eqref{eq:grcxv} as $$(8; 1^3qq1^2q1^4q1^2q).$$
  Notice that $a = 11$ and $b = 5$.  The permutation $\tau$ obtained
  by reading the indices of the $V_j$ clockwise starting with $V_1$ is
  written in one-line notation as $\{1,3,5,2,4\}$.  The values of the
  $d_j$ and $e_j$ are not illustrated in the figure, but we mention,
  for example, that $d_3 = 21$ and $e_3 = 23$.  We have also shown
  how, for some potential $C'$, that $C'(6) \in \{y:
  \rarrow{e_{\tau(5)}ye_{\tau(1)}}\}$ (as $6$ is a $1$-step for
  $C'$), but that $C'(7)$ is clearly forced to be in $V_{\tau(5)+1} =
  V_{4 + 1} = V_5$.
\end{exmp}

\begin{lemma}\label{lem:sign}
  If $\sigma\in \tpqrs$, we must have $k=\gcd(r,s,\ell)$, $r_i=r/k$
  for all $i$, and $s_i=s/k$ for all $i$.  Thus, the cycle structure
  of all elements of $\tpqrs$ is uniquely determined by $p$, $q$, $r$,
  and $s$.  Also, $\sgn(\sigma)=(-1)^{r+s+\gcd(r,s,\ell)}$.
\end{lemma}

\begin{proof}
  Take any $\sigma\in \tpqrs$. Since $\sum_{i=1}^k r_i=r$ and
  $\sum_{i=1}^k s_i=s$, Lemma \ref{lem:grinorder} implies that we must
  have $r_i=r/k$ and $s_i=s/k$ for all $i$. Then, for each $i$,
  \[ \ell_i =(r_i+s_iq)/p=\frac{(r+sq)/p}{k}=\ell/k. \]
  Note that $r_i$ and $s_i$ and (by Lemma \ref{lem:nonem}) $\ell_i$
  are all integers. By Lemma \ref{lem:gcdone}, $\gcd(r_i,s_i,\ell_i)=1$.
  Therefore
  \[ k=k\gcd(r_i,s_i,\ell_i)=\gcd(kr_i,ks_i,k\ell_i)=\gcd(r,s,\ell). \]

  The last statement of the lemma follows by noting that the sign of
  $\sigma$ is the parity of the number of letters in its domain minus
  the number of cycles in $\sigma$, including $1$-cycles. There are
  $p-r-s$ $1$-cycles, so
  \[ \sgn(\sigma)=(-1)^{p-(k+p-r-s)}=(-1)^{r+s+\gcd(r,s,\ell)}. \]
\end{proof}

We point out the fact that if $p$ is odd then the sign of $\sigma$ is
$-1$ iff $r$ and $s$ are odd.  (Note that Codenotti \& Resta
\cite[Cor. 9]{cr} determined the fact that all $\sigma\in\tpqrs$ have
the same sign when $p$ is prime.)

\section{Construction of elements in $\boldsymbol{\tpqrs}$}
\label{sec:exis}


Assume $r + sq = \ell p$ and $\gcd(r,s,\ell) = 1$.  Consider a lattice
path
\begin{equation*}
  \nu = [\nu_0 = (0,0), \nu_1, \nu_2,\dotcoms \nu_{r+s} = (r,s)]
\end{equation*}
from $(0,0)$ to $(r,s)$, where $\nu_i-\nu_{i-1}$ equals $(1,0)$ or
$(0,1)$ for $i>0$.  Associate with $\nu$ a cycle $(x;v)$ in which $v$
is an $(r+s)$-tuple in $\{1,q\}^{r+s}$ (having $r$ $1$'s and $s$
$q$'s) as follows: If $\nu_i - \nu_{i-1} = (1,0)$, then let the $i$-th
entry in $v$ be a $1$; if $\nu_i - \nu_{i-1} = (0,1)$, then set the
$i$-th entry in $v$ to be a $q$.  We refer to these cases as ``east''
and ``north'' steps, respectively.  We aim to show that if $\nu$ is
chosen appropriately, then $(x;v)$ is a well-defined element of
$\tpqrs$ for each $x\in [p]$.  It is interesting to note that our
construction of $\nu$ depends only on $r$ and $s$.

To determine $\nu$, start by setting $\nu_0 = (0,0)$ as above.
Suppose the point $\nu_i = (x_i,y_i)$ is determined.  Then set
\begin{equation}\label{eq:grcons}
  \nu_{i+1} =
  \begin{cases}
    \nu_i + (1,0), \text{ if } s x_i \leq r y_i,\\
    \nu_i + (0,1), \text{ if } s x_i > r y_i.
  \end{cases}
\end{equation}

\myfig{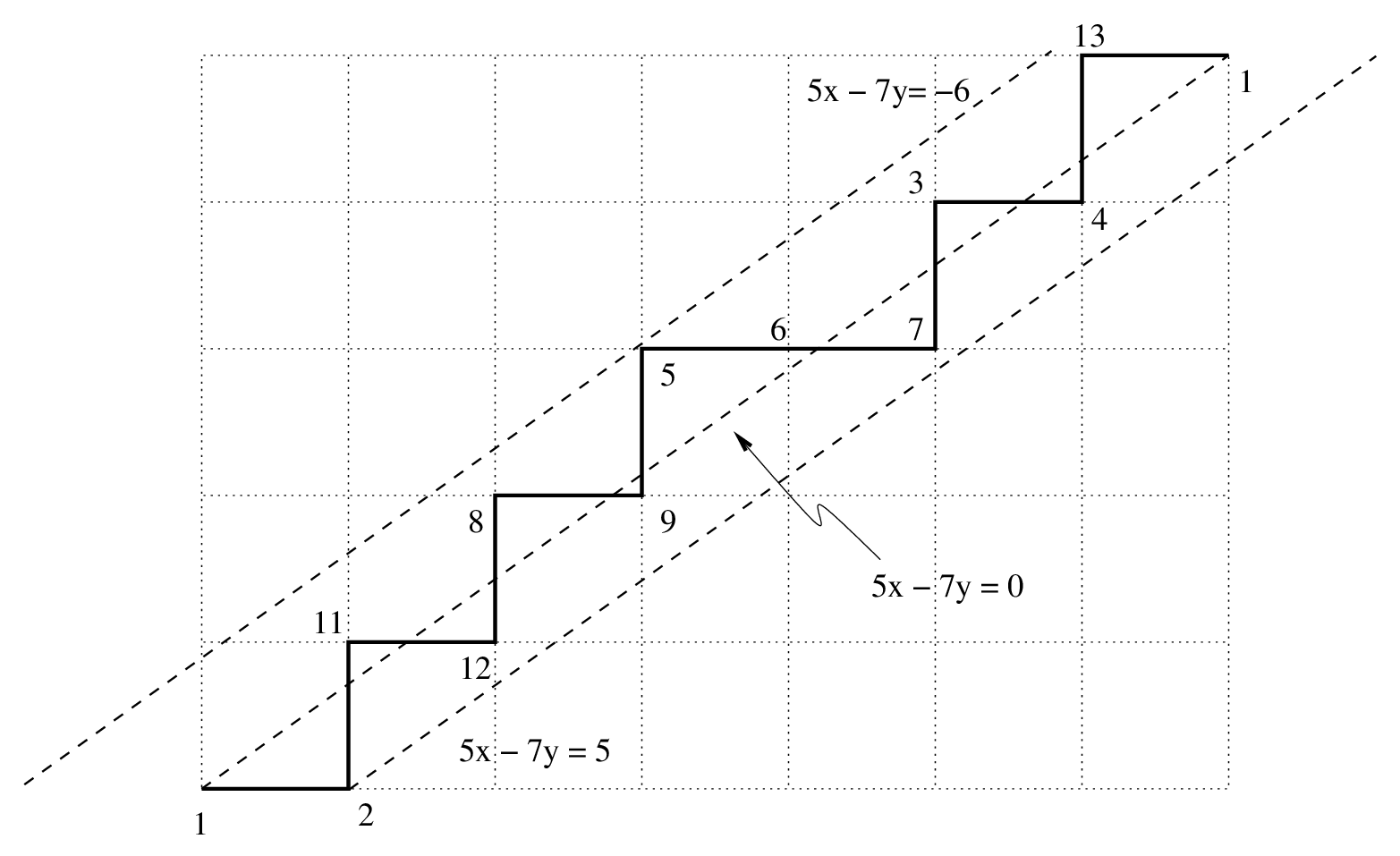}{The path $\nu$ for $p = 13$, $q = 9$, $r = 7$ and $s = 5$.}

In other words, go east if we are weakly above the line $sx-ry=0$ and
go north otherwise.  (This is effectively the Freeman approximation
used to draw diagonal lines on a computer screen.  As such, the word
$v$ can also encode the continued fraction expansion for $r/s$; see
\cite{mcilroy}.)  Figure \ref{fig:lpath} gives an example of the
construction.  In the figure, $\nu_0 = (0,0)$ is labeled by $x=1$.
Each successive $\nu_i$ is labeled by the label of $\nu_{i-1}$ plus
either $1$ or $q$ according to whether an east or north step,
respectively, separates the two vertices.  Naturally, these labels are
reduced modulo $p$.  Then the label of $\nu_i$ is precisely
$(x;v)^i(x)$.  The pair $(x;v)$ is a well-defined cycle if and only if
the only two vertices $\nu_i$ with equal labels are $\nu_0$ and
$\nu_{r+s}$.

We first bound the number of $1$-steps and $q$-steps that can appear
between any two vertices $\nu_i$ and $\nu_j$.

\begin{lemma}\label{lem:grbound}
  Determine $\nu$ as in \eqref{eq:grcons}.  Let $0\leq i , j \leq r+s$
  and write $\nu_i = (x_i,y_i)$ and
  $\nu_j = (x_j,y_j)$.  If $b = y_j - y_i$ and $a = x_j - x_i$, then
  $|as - br| \leq r + s-1$.
\end{lemma}
\begin{proof}
We claim that $-r<sx_i-ry_i\leq s$ for all points $(x_i,y_i)$ on the
path $\nu$. This is true when $i=0$, since $(x_i,y_i)=(0,0)$.
Assume the claim is true for some $i$, and consider two cases.
First, if $sx_i-ry_i\leq 0$, then $(x_{i+1},y_{i+1})=(x_i+1,y_i)$.
In this case, $sx_{i+1}-ry_{i+1}=(sx_i-ry_i)+s$, so the claim is
true for $i+1$.
Second, if $sx_i-ry_i> 0$, then $(x_{i+1},y_{i+1})=(x_i,y_i+1)$.
In this case, $sx_{i+1}-ry_{i+1}=(sx_i-ry_i)-r$, so the claim is
true for $i+1$.

Using the claim for the points $(x_i,y_i)$ and $(x_j,y_j)=(x_i+a,y_i+b)$,
we get
\[ -r+1\leq s(x_i+a)-r(y_i+b)\leq s \]
\[ -s\leq -sx_i+ry_i\leq r-1. \]
Adding gives
\[ -(r+s-1)\leq sa-rb \leq r+s-1, \]
or equivalently $|as-br|\leq r+s-1$.
\end{proof}

\begin{lemma}\label{lem:new-ineq}
  If $a,b,r,s,p$, and $q$ are integers such that $p$ divides both
  $a+bq$ and $r+sq$, then $sa-rb=0$ or $|sa-rb|\geq p$.
\end{lemma}
\begin{proof}
  Pick integers $\ell$ and $m$ such that $a+bq=pm$ and $r+sq=p\ell$.
  Then
  \[ |sa-rb|=|s(a+bq)-b(r+sq)|=|p(sm-b\ell)|. \]
  The integer $|sm-b\ell|$ is either $0$ or at least $1$, which gives
  the desired result.
\end{proof}

\begin{theorem}\label{thm:grkeq1}
  $(x;v)$ is a well-defined cycle with $r$ $1$-steps and $s$
  $q$-steps.
\end{theorem}
\begin{proof}
  $(x;v)$ has the requisite number of $1$-steps and $q$-steps by
  construction. The elements of $[p]$ moved by $(x;v)$ are those of
  the form $x+x_i+qy_i \pmod{p}$ for $0\leq i<r+s$.  We need only show
  that these $r+s$ elements are all distinct.  If this were not so,
  choose $i<j$ in the stated range with $x+x_i+qy_i\equiv
  x+x_j+qy_j\pmod{p}$. Setting $a=x_j-x_i$ and $b=y_j-y_i$ as in Lemma
  \ref{lem:grbound}, we would then have $p$ dividing $a+bq$; say, $a +
  bq = mp$. Also, by Lemma \ref{lem:new-ineq}, either $|sa-rb|=0$ or
  $|sa-rb|\geq p$.  On the other hand, Lemma \ref{lem:grbound} gives
  $|sa - rb| < r + s\leq p$.  Together, these force $sa-rb=0$.  Now,
  $b\neq 0$; otherwise $a=0$ also, contradicting the fact that
  $(x_i,y_i)\neq (x_j,y_j)$.  So we can write $r/s = a/b$ where
  $a+b<r+s$.  Let $t = \alpha/\beta \in \mathbb{Q}$ such that $r =
  at$, $s = bt$.  Pick $\alpha,\beta$ such that $\alpha,\beta \geq 1$
  and $\gcd(\alpha,\beta) = 1$.  Then
  \begin{equation}\label{eq:1}
    at + btq = \alpha\left(\frac{a}{\beta}\right) +
    \alpha\left(\frac{b}{\beta}\right)q =
    r + sq = \ell p = \alpha\left(\frac{m}{\beta}\right)p.
  \end{equation}
  Now, $\beta r = a\alpha$.  Since $\alpha$ and $\beta$ are relatively
  prime, we conclude that $\beta$ divides $a$.  Similarly, $\beta$
  divides both $b$ and $m$.  So from \eqref{eq:1}, $\alpha$ divides\
  $r,s$, and $\ell$.  As $a < r$, we must have $\alpha > \beta \geq
  1$.  This yields a contradiction with our requirement that
  $\gcd(r,s,\ell) = 1$.
\end{proof}

We now relax the assumption that $\gcd(r,s,\ell) = 1$.  Indeed,
suppose this $\gcd$ is $k > 1$.

Consider $(x;v)$ where $v$ is determined by the lattice path $\nu$
from $(0,0)$ to $(r/k,s/k)$ constructed in \eqref{eq:grcons}.  Theorem
\ref{thm:grkeq1} assures us that $(x;v)$ is a valid cycle.

\begin{theorem}\label{thm:nonempty}
  Let $k=\gcd(r,s,\ell)$ and $\nu$ be as above and write $C_j$ for
  $(1+(j-1)(q-1);v)$.  Then 
  \[\sigma = C_1C_2\cdots C_k\]
  is well-defined element of $\tpqrs$.
\end{theorem}
\begin{proof}
  We already know that each cycle $C_j$ is well-defined; it suffices
  to check that these cycles are disjoint. The set
\[ \{1+(j-1)(q-1)+x_i+qy_i \pmod{p}: 0\leq i<r/k+s/k,\,1\leq j\leq k\}. \]
consists of those elements moved by $C_j$.  Suppose two such elements
are equal mod $p$, say
\[ 1+(j_1-1)(q-1)+x_{i_1}+qy_{i_1}=
   1+(j_2-1)(q-1)+x_{i_2}+qy_{i_2}+pM. \]
We must show that $i_1=i_2$ and $j_1=j_2$. Choose labels so that $j_1\geq
j_2$.
Set $a=x_{i_2}-x_{i_1}$, $b=y_{i_2}-y_{i_1}$, $\rpri =r/k$, $\spri =s/k$,
$\ell'=\ell/k$, and $j=j_1-j_2$. We then have $0\leq j\leq k-1$ and
\[  j(q-1)=a+qb+pM. \]
Set $A=a+j$ and $B=b-j$. Then $p(-M)=A+qB$, so that $p$ divides
$A+qB$. Since $p$ also divides $\rpri +s'q$, Lemma \ref{lem:new-ineq}
says that $\spri A-\rpri B=0$ or $|\spri A-\rpri B|\geq p$.

Assume the second alternative occurs. Then
\[ |\spri a-\rpri b+j(\spri +\rpri )|\geq p. \]
Now Lemma \ref{lem:grbound} gives
\[ |\spri a-\rpri b|<\rpri +\spri . \]
Hence,
\[ j(\spri +\rpri )\geq |\spri a-\rpri b+j(\spri +\rpri )|-|\spri a-\rpri b| > p-(\rpri +\spri ). \]
This gives $j>\frac{p}{\rpri +\spri }-1$.
But $r+s\leq p$, so that $\rpri +\spri =r/k+s/k\leq p/k$, which implies
$k\leq\frac{p}{\rpri +\spri }$. We deduce that $j>k-1$, contradicting
the fact that $0\leq j\leq k-1$.

We must therefore have $\spri A-\rpri B=0$, or $\spri a-\rpri b=-j(\spri +\rpri )$.
It is still true that $|\spri a-\rpri b|<\rpri +\spri $, so we see that
\[ |j(\rpri +\spri )|<\rpri +\spri . \]
Since $j$ is an integer and $\rpri +\spri >0$, we must have $j=0$ and $j_1=j_2$.
Then $\spri a-\rpri b=0$ as well. If $a=b=0$, then $i_1=i_2$ and we are done.
Otherwise, both $a$ and $b$ are nonzero and we get $\rpri /\spri =a/b$
with $a+b<\rpri +\spri $. This contradicts $\gcd(\rpri ,\spri ,\ell')=1$, just as in
the proof of Theorem \ref{thm:grkeq1}.
\end{proof}

\myfig{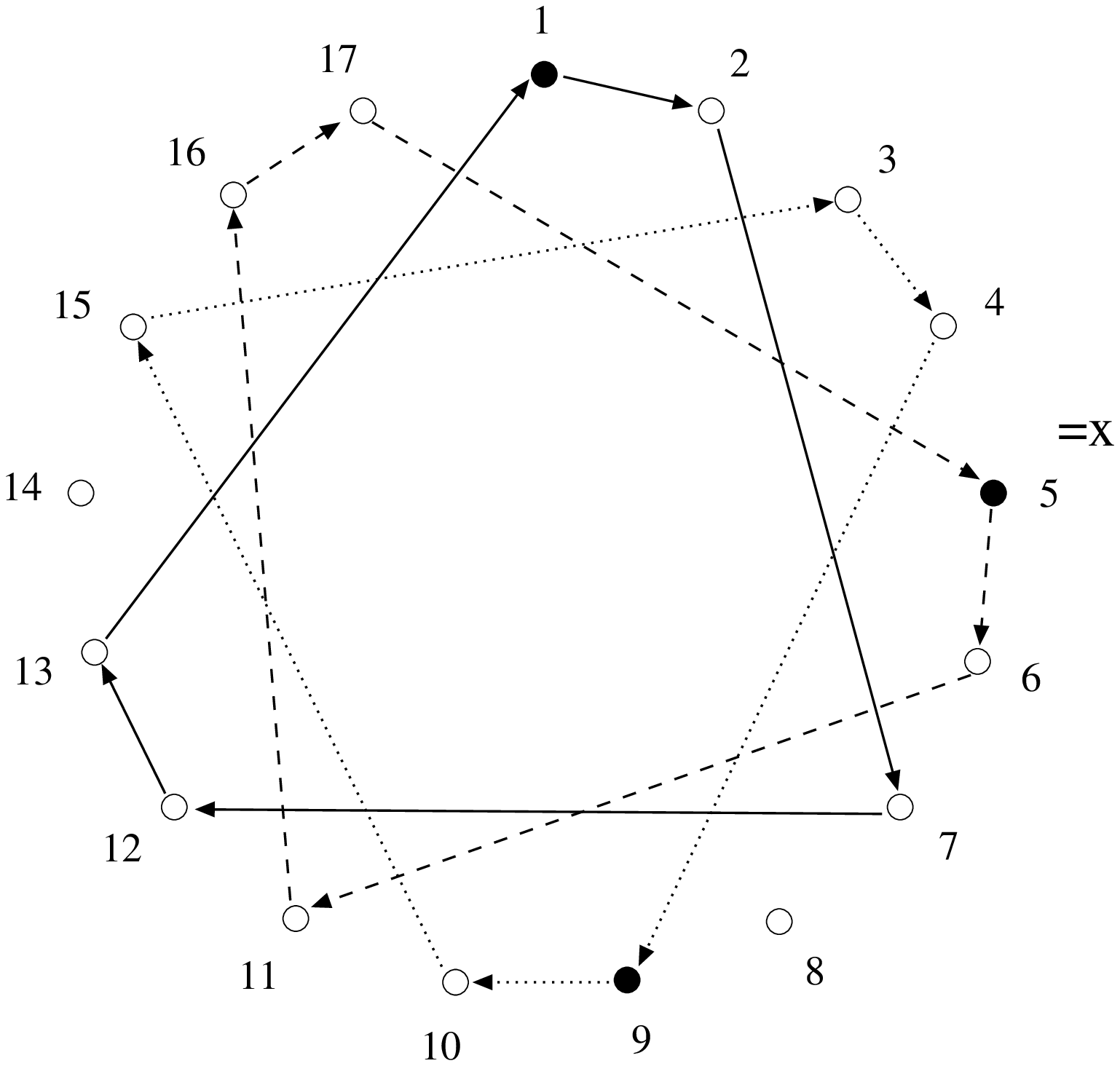}{Illustration for Example \ref{ex:bigk}}
\begin{exmp}\label{ex:bigk}
  We illustrate the case of $p = 17$, $q = 5$, $r = 6$ and $s = 9$.
  $r + sq = 6 + 9\cdot 5 = 51 = 3\cdot 17$, so $\ell = 3$ and $k =
  \gcd(r,s,\ell) = 3$.  Shown are $C_1 = (1; v)$ (solid), $C_2 = (1 +
  4; v)$ (dashed) and $C_3 = (1 + 2\cdot 4; v)$ (dotted).
\end{exmp}

\begin{theorem}
The coefficient $a(r,s)$ in the circulant determinant is
zero if $p$ does not divide $r+qs$. Otherwise, this coefficient
is nonzero with sign $(-1)^{\gcd(r,s,(r+qs)/p)}$.
\end{theorem}
\begin{proof}
  Immediate from all the preceding results.
\end{proof}

\section{The largest coefficient}
\label{sec:maxc}
We have identified the coefficients of the monomials in $\Phi_{p,q}$ as the
numbers of permutations in certain classes. In this section we will obtain
two-sided bounds on the size of the largest coefficient.

Consider the \textit{permanent} of the circulant matrix
\[D_{p,q}(x,y)= \circu(1,x,0,\dots 0,y,0,\dots,0),\]
in which the $y$ is the $(q+1)^{\mathrm{st}}$ entry. Since all of the
permutations that contribute to a given monomial in the determinant
\[\Phi_{p,q}(x,y)=\det{(\circu(1,-x,0,\dots 0,-y,0,\dots,0))},\]
have the same sign, it follows that if
\[\Phi_{p,q}(x,y)=\sum_{r,s}a_{p,q}(r,s)x^ry^s,\]
then
\[D_{p,q}(x,y)=\sum_{r,s}|a_{p,q}(r,s)|x^ry^s.\]
Thus $D_{p,q}(1,1)$ is the sum of the absolute values of the coefficients
$a_{p,q}(r,s)$. Let $M(p,q)=\max_{r,s}|a_{p,q}(r,s)|$. Then we have
\[  \frac{D_{p,q}(1,1)}{N(p,q)}    \le M(p,q)\le D_{p,q}(1,1),\]
in which $N(p,q)$ is the number of distinct monomials that appear.

We now obtain two-sided estimates for $D_{p,q}(1,1)$. This is the permanent
of a circulant matrix that has three cyclic diagonals of 1's and whose other
entries are 0's.

For the upper bound we have the following theorem of Br\`egman-Minc
\cite{bre, min,sch}.
\begin{theorem}
[Br\`egman, Minc] Let $A$ be an $n\times n$ 0-1 matrix with $r_i$ 1's in row
$i$, for each $i=1,2,\dots ,n$. Then the permanent of $A$ satisfies
\[\per(A)\le \prod_{i=1}^n(r_i!)^{1/r_i},\]
and the sign of equality holds iff $A$ consists of a sequence of $r_i\times
r_i$ blocks of 1's on the main diagonal, with all other entries being 0's.
\end{theorem}
If we apply this theorem to $D_{p,q}(1,1)$ we find that
\[M(p,q)\le 6^{p/3}= (1.817..)^p.\]

For the lower bound we have the theorem of Egorychev \cite{ego} and van der
Waerden.
\begin{theorem}[van der Waerden, Egorychev] Let $A$ be an $n\times n$ matrix
whose entries are nonnegative and sum to 1 in every row and column. Then
$\per(A)\ge n!/n^n$, with equality iff $A$ is the matrix whose entries are
all equal to $1/n$.
\end{theorem}
We apply this theorem to $\circu(1,1,0,\dotcoms 0,1,0,\dots,0)$/3. The
result is that
\[D_{p,q}(1,1)\ge \frac{3^pp!}{p^p}\sim \left(\frac{3}{e}\right)^p\sqrt{2\pi
e}.\]

Finally since $N(p,q)$, the number of monomials that appear, is at most
$p^2$, we have proved the following.
\begin{theorem}
Fix $q$. Then the maximum absolute value of the coefficients in the
polynomial $\Phi_{p,q}(x,y)$ satisfies
\[1.1036..=\frac{3}{e}\le \liminf_{p\to\infty}M(p,q)^{1/p}\le
\limsup_{p\to\infty}M(p,q)^{1/p}\le 6^{1/3}=1.817...\]
In particular, the largest coefficient grows exponentially with $p$.
\end{theorem}

\bigskip

\textbf{Acknowledgment:} The authors thank John D'Angelo for useful
discussions involving this problem.

\newpage


\end{document}